# CRITICAL PERCOLATION OF VIRTUALLY FREE GROUPS AND OTHER TREE-LIKE GRAPHS

By Iva Špakulová

*Vanderbilt University*

This article presents a method for finding the critical probability $p_c$ for the Bernoulli bond percolation on graphs with the so-called tree-like structure. Such a graph can be decomposed into a tree of pieces, each of which has finitely many isomorphism classes. This class of graphs includes the Cayley graphs of amalgamated products, HNN extensions or general groups acting on trees. It also includes all transitive graphs with more than one end.

The idea of the method is to find a multi-type Galton–Watson branching process (with a parameter $p$) which has finite expected population size if and only if the expected percolation cluster size is finite. This provides sufficient information about $p_c$. In particular, if the pairwise intersections of pieces are finite, then $p_c$ is the smallest positive $p$ such that $\det(M-1) = 0$, where $M$ is the first-moment matrix of the branching process. If the pieces of the tree-like structure are finite, then $p_c$ is an algebraic number and we give an algorithm computing $p_c$ as a root of some algebraic function.

We show that any Cayley graph of a virtually free group (i.e., a group acting on a tree with finite vertex stabilizers) with respect to any finite generating set has a tree-like structure with finite pieces. In particular, we show how to compute $p_c$ for the Cayley graph of a free group with respect to any finite generating set.

**1. Introduction.** We will use the notation $\mathcal{G} = (V, E)$ for a graph with the vertex set $V$ and the edge set $E$. All graphs are assumed to be locally finite (vertices have finite degrees) and transitive (for any two vertices, there exists an automorphism of $\mathcal{G}$ mapping one to the other). We fix one vertex of the graph and call it the *origin*.

For every $p \in (0, 1)$, the *Bernoulli bond percolation* on $\mathcal{G}$ is a product probability measure $P_p$ on the space $\Omega = \{0, 1\}^E$ of subsets of the edge set $E$. The product measure is defined via $P_p(\omega(e) = 1) = p$ for all $e \in E$. The $\sigma$-algebra of $P_p$-measurable sets does not depend on $p$. We denote the $\sigma$-algebra by $\Sigma$ and the expected value by $E_p$.

For any realization $\omega \in \Omega$, open edges form a random subgraph of $\mathcal{G}$. The *percolation function* is defined to be the probability that the origin is contained in an infinite cluster. The behavior of the percolation model depends strongly on the value of probability $p$. There is a critical value $p_c$ of the probability $p$ such that, for $0 \leq p < p_c$, all clusters are finite, and if $p_c < p \leq 1$, there is an infinite cluster $P_p$-almost surely (see Grimmett [4]).

Explicit values of $p_c$ have been found only for some special cases. In particular, for lattices in $\mathbb{R}^2$, the value of $p_c$ is obtained using dual graphs (for $\mathbb{R}^d$ with $d \geq 3$, the values of $p_c$ are not known). Kesten [5] proved, for the square lattice, that $p_c = 1/2$. For the triangular lattice that $p_c = 2\sin(\pi/18)$; for the hexagonal lattice, that $p_c = 1 - 2\sin(\pi/18)$ (see Grimmett [4]). Ziff and Scullard [11] found $p_c$ for a larger class of lattices in $\mathbb{R}^2$ (they considered graphs that can be decomposed into certain self-dual arrangements). The critical probability value is also known for trees ($p_c = 1/\text{branching number}$) and virtually cyclic groups, where $p_c = 1$. The critical probability can also be found in the case of a free product of finite transitive graphs [6]. As far as we know, these are the only known graphs where $p_c$ for the bond percolation has been computed exactly.

If we change the generating set of a group, the graph changes dramatically, and there have been no examples of groups (except virtually cyclic ones) where $p_c$ was known for all generating sets. It is expected that many properties of percolation (behavior at the critical value) are invariant with respect to changing the generating set, so it is useful to know how $p_c$ depends on the generating set of a group.

In this article, we study the critical probability $p_c$ for a class of graphs that admit the so-called *tree-like structure*. Roughly speaking, such a graph can be decomposed into a rooted tree of edge-disjoint pieces that intersect at the so-called border sets, which are cut sets of the graph. We always assume that there are finitely many isomorphism classes of pieces. This class of graphs includes, for example, all transitive graphs with more than one end, Cayley graphs of amalgamated products and HNN extensions [e.g., $SL(2,\mathbb{Z}) = \mathbb{Z}_4 *_{\mathbb{Z}_2} \mathbb{Z}_6$]. In the case of an amalgamated product, the pieces correspond to the Cayley graphs of factor groups and the border sets consist of vertices from cosets of the amalgamated subgroup. The precise definition of the tree-like structure is contained in Section 2 and several examples, including the grandparent tree, can be found in Section 5.

Consider a realization of the percolation on a graph $\mathcal{G}$ with a tree-like structure. For each piece $P_i$ that is not the root, let $B_i$ be the border set that is the intersection of $P_i$ with its parent in the tree.

Since $B_i$ is a cut set, we can define an equivalence relation on it by saying that two vertices of $B_i$ are equivalent if they are connected by an open path inside the union of all pieces that are not descendants of $P_i$. The data consisting of the equivalence relation and the distinguished equivalence class connected to the origin is called the *color* of the piece $P_i$. If no vertex of the border set is connected to the origin, then we say that this piece has the color white. Since the $p_c$ of the whole graph does not exceed the $p_c$ of any subgraph, we can assume that the probability $p$ is smaller than the minimum of all $p_c$'s of pieces. It is then easy to see that the percolation process dies if and only if the tree of nonwhite pieces is finite. The fact that the pieces form a tree suggests the use of a branching process with individuals corresponding to the colored pieces, such that the distribution of children is induced by the percolation process. Unfortunately, the color of a piece depends not only on the color of its parent (as required for a branching process), but also on the colors of its siblings and their descendants.

Nevertheless, we define a different distribution on the colors of the children and we obtain a branching process that has finite population size if and only if the expected size of the percolation cluster is finite. The next statement is the main result of the paper.

THEOREM 1.1.   *Assume that the graph $\mathcal{G}$ has a tree-like structure:*

  (i) *For a percolation with parameter $p$, there exists a branching process on the tree of pieces such that the expected size of its population is finite if and only if the expected size of the percolation cluster at the origin is finite.*

  (ii) *If all of the border sets are finite, then the branching process has finitely many types and the first moment matrix is of finite size. In this case, $p_c$ is the smallest value of $p$ such that the spectral radius of the first moment matrix is 1.*

  (iii) *If, in addition, the pieces are finite, then the entries of the first moment matrix are algebraic functions in $p$. Therefore, $p_c$ is algebraic.*

*There exists an algorithm that, given the pieces and their border sets, computes a finite extension $K$ of the field $\mathbb{Q}(p)$ and an algebraic function $f$ in $K$ such that $p_c$ is the smallest positive root of $f$.*

The theorem is proved in Section 4.

We list here two corollaries of Theorem 1.1. The first has already been proven in [6], by the author, using different methods. The second corollary answers a question of M. Sapir (personal communication) about the special linear group $SL(2,\mathbb{Z})$, a simple example of a group that is not free product. Both corollaries are proved in Section 5.

By the expected subcritical cluster size $\chi_i(p)$, we mean $E_p(|C|)$, where $C$ is the cluster containing the origin and $p < p_c$.

COROLLARY 1.2. *Let $\mathcal{G}$ be a free product of (transitive) graphs. Denote by $\chi_i(p)$ the expected (subcritical) cluster size in the ith factor graph. The critical probability $p_c$ of $\mathcal{G}$ is the infimum of positive solutions of*

$$\sum_{j=1}^{n} \prod_{i=1, i\neq j}^{n} \chi_i(p) - (n-1) \prod_{i=1}^{n} \chi_i(p) = 0.$$

Note that in the case of free products (with respect to natural set of generators), the border sets consist of one vertex and so the branching process has just one type of individual (for more details, see Section 5.1). It is, in fact, also a special case of Theorem 1.4 below.

We explicitly compute the critical probabilities of several Cayley graphs. In particular, we prove the following.

COROLLARY 1.3. *The critical probability $p_c$ of $SL(2,\mathbb{Z})$ given by presentation $\langle a,b|a^4,b^6,a^2b^{-3}\rangle$ is an algebraic number that is equal to $0.4291140496\ldots$.*

The general case of Cayley graphs of amalgamated products and HNN extensions is covered by the following theorem. Consider a group $G$ acting on a simplicial tree $T$. The standard generating set of the group $G$ is any generating set consisting of elements in the vertex stabilizers and free letters. It follows from the structure theorems of Bass–Serre theory (see, e.g., [2] or [10]) that the group $G$ is the fundamental group of a graph of groups (see Section 6.3).

THEOREM 1.4. *Let $G$ be a group acting on a simplicial tree and let $\mathcal{G}$ be its Cayley graph with respect to a standard generating set. Then, $\mathcal{G}$ has a tree-like structure whose pieces correspond to the Cayley graphs of the vertex stabilizers and whose border sets correspond to the edge stabilizers.*

In the last section of the paper, we prove that every transitive graph with more than one end has a tree-like structure with finite border sets. This can be applied to, say, Cayley graphs of free groups with arbitrary finite generating sets. Moreover, the pieces obtained from the general construction are finite in this case, and are explicitly described in Section 7.2.

THEOREM 1.5. *Let $G$ be a virtually free group, that is, it acts on a simplicial tree with finite vertex stabilizers. Its Cayley graph with respect to any finite generating set then has a tree-like structure with finite border sets and finite pieces. Given a finite generating set, the pieces of the tree-like structure are algorithmically constructed.*

*Therefore, the $p_c$ is an algebraic number and one can use the algorithm from Theorem 1.1(*iii*) to compute $p_c$, given any finite generating set.*

This theorem is proved in Section 6.5.

This provides the first example of a class of Cayley graphs closed under quasi-isometry where we can algorithmically find the value of $p_c$ for every graph in the class (besides the graphs with 0 and 2 ends, where $p_c = 1$).

**2. Graphs with tree-like structure.** In what follows, all graphs are supposed to be transitive, locally finite, connected and infinite (e.g., Cayley graphs of finitely generated groups). By $V(\mathcal{G})$ [resp., $E(\mathcal{G})$], we denote the set of vertices (resp., edges) of a graph $\mathcal{G}$.

Here, we present the formal definition of the tree-like structure and several basic properties. The reader may wish to review the examples presented in Section 5 to obtain a more intuitive picture.

2.1. *The definition.* We are now going to define a tree-like structure of a graph.

Note that a rooted tree of elements defines a partial order on these elements. In this partial order, an element is smaller than another if and only if it is its descendant. The root is a maximal element in such a case.

DEFINITION 2.1. A *tree-like structure* on a (transitive, locally finite, connected, infinite) graph $\mathcal{G}$ is a triple $(\mathbf{P}, J, \gamma)$, where $\mathbf{P}$ consists of pairs of nonempty subgraphs $(P_i, B_i)$ of $\mathcal{G}$, $i \in I$ ($P_i$ are called the *pieces*, $B_i$ are called the *border sets*), $J$ is a finite subset of $I$, $\gamma$ is a *model map* from $I$ to $J$ and the following conditions are satisfied.

(1) For every $i \in I$, $P_i$ is a subgraph of $\mathcal{G}$ and $B_i \subseteq V(P_i)$.
(2) For every $i \neq j$, $E(P_i) \cap E(P_j) = \varnothing$ and $\bigcup_{i \in I} E(P_i) = E(\mathcal{G})$.
(3) There is a partial order on the pieces with maximal element $P_0$ such that its graphical representation is a tree with a root $P_0$; moreover, if $P_i$ is a child of $P_j$, then $P_i \cap P_j = B_i$.

[Denote by $U(P_i)$ the union of the pieces in the descendant subtree of $P_i$ ($P_i$ included) for all $i \in I$.]

(4) For every $i \neq j$, if $P_i \nsubseteq U(P_j)$, then $P_i \cap U(P_j) \subseteq B_j$.
(5) For every $i \in I$, there is an isomorphism between $U(P_i)$ and $U(P_{\gamma(i)})$ taking pieces to pieces and border sets to border sets, respecting the order on the pieces.

We say that there is a finite number of isomorphism classes of pairs $(P_i, B_i)$. For $j \in J$, the pieces $P_j$ are called the *model pieces*.

If there is a piece $P_i$ with no edges, then it consists only of vertices and we can remove such a piece (and add the vertices to the parent piece) and change the tree-like structure accordingly. Therefore, we will assume that each $P_i$ contains at least one edge. Clearly, the tree-like structure can be degenerated in the sense that the whole graph is just one piece or the number of pieces is finite. In what follows, we will always assume that it is nondegenerate, that is, that the number of pieces is infinite. In all cases considered in this paper, the pieces have at most finitely many components (this simplifies the computations).

For every $i \in I$, denote by $\Lambda_i \subset I$ the set of indices of children of $P_i$ [given by part (3) of Definition 2.1]. We will now present several basic properties of the tree-like structure.

LEMMA 2.2. *For every $i$, the set $U(P_i)$ is covered by $P_i$ and the collection of $U(P_\lambda)$ for $\lambda \in \Lambda_i$. Moreover, pairwise intersections of $U(P_\lambda)$'s for $\lambda \in \Lambda_i$ consist only of vertices of $P_i$.*

PROOF. By property (2), each edge is included in exactly one piece and therefore, if $P_k \subset U(P_{\lambda_1})$, then $P_k \not\subset U(P_{\lambda_2})$, provided $\lambda_2 \neq \lambda_1$. Therefore, by property (4), the intersection $U(P_{\lambda_1}) \cap U(P_{\lambda_2})$ is included in $B_{\lambda_1} \cap B_{\lambda_2}$. Since each $B_\lambda \subset P_i$, by property (3), the claim follows. □

LEMMA 2.3. *Every border set (except possibly $B_0$, the border set of the root) is a vertex cut set of graph $\mathcal{G}$. In particular, let a subgraph $U(P_i)^c$ of $\mathcal{G}$ be induced by the edges $E(\mathcal{G}) \setminus E(U(P_i))$. Then, $U(P_i) \cup U(P_i)^c = \mathcal{G}$ and $U(P_i) \cap U(P_i)^c = B_i$.*

PROOF. If the complement $U(P_i)^c$ satisfies $U(P_i) \cap U(P_i)^c = B_i$, then $B_i$ is indeed a cut set. By property (4), for any $P_k \subset U(P_i)^c$, we have $P_k \cap U(P_i) \subset B_i$, so the claim follows. □

2.2. *Changing the generating set.* The following lemmas apply to the Cayley graphs of groups. In that case, the group itself acts on its Cayley graph (by multiplication from the right) and so it is a subgroup of the automorphism group $Aut(\mathcal{G})$. If the isomorphisms of part (5) of Definition 2.1 are (almost) in the group $G$, then we can make the generating set smaller or slightly larger and still obtain a tree-like structure.

LEMMA 2.4. *Let $S$ be a finite generating set of a group $G$. Suppose that there is a tree-like structure of the Cayley graph $\mathcal{G}$ of $G$ with respect to $S$, that each of the isomorphisms from part (5) of Definition 2.1 can be extended to the whole of $\mathcal{G}$ and that these extensions form a subgroup $H$ of $Aut(\mathcal{G})$. If $G \cap H$ has finite index in $H$, then, for any generating subset $S' \subset S$, the Cayley graph $\mathcal{G}'$ of $G$ with respect to $S'$ has a tree-like structure with pieces containing the same vertices (and fewer edges).*

PROOF. If we restrict the isomorphisms in part (5) of Definition 2.1 from $H$ to $G \cap H$, then we will increase the number of isomorphism classes by a finite amount because $G \cap H$ has finite index in $H$. Thus, we can assume that $H \subset G$.

The graph $\mathcal{G}'$ is obtained from $\mathcal{G}$ by removing edges with labels in $S \setminus S'$. The pieces $P_i'$ of the tree-like structure on $\mathcal{G}'$ are obtained in the same way (by simply removing edges with labels in $S \setminus S'$). Let us verify the conditions of Definition 2.1. Every edge belongs to exactly one piece, thus property (2) follows. Conditions (1), (3) and (4) depend only on vertices of pieces. Since vertices do not change, these conditions are satisfied.

It remains to show property (5). Since the isomorphisms between $U(P_i)$'s in $\mathcal{G}$ are given by elements in $G$, they are also isomorphisms of $U(P_i')$'s in $\mathcal{G}'$. There are thus finitely many isomorphism classes between pieces of $\mathcal{G}'$. □

LEMMA 2.5. *Let $S$, $S'$ be finite generating sets of a group $G$. Suppose that there is a tree-like structure with pieces $P_i$ of the Cayley graph $\mathcal{G}$ of $G$ with respect to $S$ and that the isomorphisms from part (5) of Definition 2.1 form a subgroup $H$ of $Aut(\mathcal{G})$. Assume that $G \cap H$ has finite index in $H$ and that the following condition is satisfied:*

(∗) *for all $x \in S' \setminus S$ and $g \in G$, there is a piece $P_i$ containing the two vertices labeled $g$ and $gx$.*

*The Cayley graph $\mathcal{G}'$ of $G$ with respect to $S'$ then has a tree-like structure with pieces $P_i'$, where $P_i'$ contains the same vertices as $P_i$.*

PROOF. We can again assume that $H \subset G$. If we prove the claim for $S'' = S \cup S'$, then, for $S'$, it follows from Lemma 2.4. Therefore, we can assume that $S \subset S'$.

The graph $\mathcal{G}'$ is obtained from $\mathcal{G}$ by adding edges labeled by elements in $S' \setminus S$. In the tree-like structure, we add each new edge to the oldest (i.e., the maximal in the partial order) piece $P_i$ containing both of its endpoints. Such a piece exists by assumption (∗). Suppose that there are two such oldest pieces, $P_i$ and $P_j$, containing endpoints of an edge $e$. Then, neither of them is a descendant of the other and their intersection is in $B_i \cup B_j$, by condition (4). Therefore, $e$ should be added to their predecessor that is older than both of them. Therefore, there is a unique such oldest piece for each edge. This implies property (2). Conditions (1), (3) and (4) now depend only on vertices of pieces that have not been changed. Therefore, these conditions hold as in the tree-like structure on $\mathcal{G}$.

It remains to prove property (5). Consider an isomorphism from property (5) of the graph $\mathcal{G}$. Let $h \in H \subset G$ such that $U(P_i) = hU(P_j)$ and $h$ takes pieces to pieces and border sets to border sets (we can assume that $i, j \neq 0$).

We will show that the descendant trees of the modified pieces $P'_i$ and $P'_j$ in $\mathcal{G}'$ can be mapped by the same isomorphism $h$ and that it also takes the modified pieces to pieces and border sets to border sets. The vertices of pieces have not been changed, so a difficulty can arise only for edges. By contradiction, assume that $U(P'_i) \neq hU(P'_j)$ that means that there is a descendant $P'_k$ of $P'_i$ and $P'_l$ of $P'_j$ such that $P_k = hP_l$ but $P'_k \neq hP'_l$. There exists an edge $e$ in $P'_k$ that is missing in $hP'_l$ (or vice versa). If $e \notin hP'_l$, then there is a older (in the ordering) piece than $P'_l$ containing both endpoints of $h^{-1}e$. Therefore, these endpoints are in the border set $B'_l$. Consequently, the endpoints of $e$ are in $B'_k$ and so the edge must be included in some ancestor of $P'_k$, a contradiction. $\square$

Lemma 2.5 can be applied only if the pieces are large enough to satisfy the condition $(*)$. The following lemma allows us to enlarge pieces. The new piece $P'_i$ is a union of the children of the piece $P_i$. If we started with finite pieces, then the constructed pieces are again finite. Moreover, the tree of pieces (partial order) remains unchanged.

LEMMA 2.6. *Assume that graph $\mathcal{G}$ admits a tree-like structure with pieces $P_i$ and border sets $B_i$. Define new pieces $P'_i$ and border set $B'_i$ as follows:*

$$P'_0 = P_0 \cup \bigcup_{k \in \Lambda_0} P_k, \qquad B'_0 = B_0,$$

$$P'_i = V(P_i) \cup \bigcup_{k \in \Lambda_i} P_k, \qquad B'_i = V(P_i) \qquad \text{for } i \neq 0.$$

*Define the partial order on the pieces $P'_i$ to be the same as on the original pieces $P_i$. The modified pieces $P'_i$ and border sets $B'_i$ then form a new tree-like structure of $\mathcal{G}$.*

PROOF. We need to verify the conditions of Definition 2.1:

(1)–(2) These clearly follow from the definitions of $P'_i$ and $B'_i$.
  (3) Let $P'_j$ be a child of $P'_i$. Then, $P'_i \cap P'_j$ contains vertices of $P_j$, that is, $B'_j$. If $k \in \Lambda_j$, then, by property (4) of the original tree-like structure, $U(P_k) \cap P'_i \subset B_k \subset V(P_j) = B'_j$. Therefore, $P'_i \cap P'_j = B'_j$.
  (4) Assume that $P'_i$ is not a descendant of $P'_j$ and that $i \neq j$. Then, $P_i$ is not a descendant of $P_j$ and no child of $P_i$ is a descendant of $P_j$ (it can be $P_j$ itself). Therefore, using property (4) of the original tree-like structure, we see that $P'_i \cap U(P_j) \subset P_j$. Since the intersection contains only vertices [by (2)], we can write

$$P'_i \cap U(P'_j) = V\left(\left(P_i \cup \bigcup_{k \in \Lambda_i} P_k\right) \cap U(P_j)\right) \subset V(P_j) = B'_j.$$

(5) Assume that there exists an isomorphism, from property (5), between $U(P_i)$ and $U(P_j)$ such that the pieces and border sets are respected. The modification of the piece $P_i$ into $P'_i$ uses the children of $P_i$ that are preserved by the isomorphism. Thus, the modified pieces and their border sets are also preserved by this isomorphism between $U(P'_i)$ and $U(P'_j)$. There are thus finitely many isomorphism classes. $\square$

**3. Partitions of the border sets.** Consider a graph $\mathcal{G}$ with a tree-like structure as above and take a realization $\omega \in \Omega$ of a percolation process on $\mathcal{G}$. Some pairs of vertices of the border set $B_0$ can be connected by open paths in $\mathcal{G}$. In this way, the realization determines a partition of the set $B_0$: two vertices are in the same class of the partition if they are connected by an open path in $\mathcal{G}$. The percolation process induces a probability measure on the set of all partitions of $B_0$. Similarly, we obtain a partition of each $B_j$ by looking at the open paths in the subgraph $U(P_j)$. We call this a *descendant partition* of the border set $B_j$. We say that a partition of a set $B$ is *induced by percolation* on $H$ if $H$ is a subgraph of $\mathcal{G}$, $B$ is a subset of vertices of $H$ and two vertices of $B$ are in the same class of the partition if and only if they are connected by an open path in $H$.

In this section, we will use the decomposition into pieces to find the measure on the set of descendant partitions of border sets using recurrent relations.

Let $Z^{(i)}$ be the set of all partitions of the border set $B_i$ that can be induced by the percolation on $U(P-i)$. Let $q_i : \Omega \to Z^{(i)}$ be a map assigning to each realization $\omega \in \Omega$ the partition on the border set $B_i$. We consider only partitions induced by percolation on $U(P-i)$ and thus the map $q_i$ is surjective. The measure on $Z^{(i)}$ is a pullback of the percolation measure $P_p$ by $q_i$. Thus, $A \subset Z^{(i)}$ is measurable if and only if $q_i^{-1}(A) \in \Sigma$. For simplicity, we use the same notation $P_p$ for the probability measure on partitions and denote the $\sigma$-algebra of measurable sets of partitions by $\mathcal{Z}^{(i)}$.

For every finite subset $F$ of $B_i$ and every partition $z^{(i)}$ of $B_i$, denote by $z^{(i)}(F)$ the set of all partitions that coincide with $z^{(i)}$ on $F$. These sets are obviously in $\mathcal{Z}^{(i)}$.

If the subgraphs $U(P_i)$ and $U(P_j)$ are isometric [by an isomorphism from part (5) of the definition], then the $\sigma$-algebras of partitions $Z^{(i)}$ and $Z^{(j)}$ are isomorphic and the measures induced by the same percolation process are preserved.

LEMMA 3.1. *A descendant partition of $B_i$ is determined by the state of the edges in $P_i$ and by the descendant partitions of border sets $B_\lambda, \lambda \in \Lambda_i$ (which are children of $P_i$).*

PROOF. Indeed, using property (2) of Definition 2.1 and Lemma 2.2, we see that every edge of an open path connecting two vertices of $B_i$ in $U(P_i)$ is in exactly one of the following graphs: $P_i$ or $U(P_\lambda), \lambda \in \Lambda_i$. We can split the path into several segments, each of them being in some $U(P_\lambda)$ or $P_i$. A segment in $U(P_\lambda)$ starts at some point of $B_\lambda$ and ends at some other point of $B_\lambda$. Therefore, the existence of the open segment is determined by the partition on $B_\lambda$. Thus, in order to decide whether there is an open path connecting certain vertices in $B_i$, it is sufficient to know the states of edges in $P_i$ and the partitions of $B_\lambda, \lambda \in \Lambda_i$. □

For every measurable set $A \in \mathcal{Z}^{(i)}$ and every vector $\zeta \in \prod_{\lambda \in \Lambda_i} Z^{(\lambda)}$ of partitions on $B_\lambda$'s, $\lambda \in \Lambda_i$, we denote the conditional probability of $A$ provided $\zeta$ by $\mathrm{P}_p(A|\zeta)$. Denote by $\mu_p^{(i)}$ the product measure on $\prod_{\lambda \in \Lambda_i} Z^{(\lambda)}$ such that the measure on each $Z^{(\lambda)}$ is given by $P_p$.

We can then write

$$\mathrm{P}_p(A) = \int_{\zeta \in \prod_{\lambda \in \Lambda_i} Z^{(\lambda)}} \mathrm{P}_p(A|\zeta) \, d\mu_p^{(i)}(\zeta). \tag{1}$$

Note that since for pieces of the same isomorphism class, we have isomorphic $\sigma$-algebras of the set of partitions and measures on them, we can consider only one such equation for each model border set (of index in the finite set $J$ from the definition).

This defines an operator on the space of all measures on the direct product of $Z^{(i)}$'s. More precisely, let $X_j$ be the space of the probability measures on $Z^{(j)}$ with the $\sigma$-algebra $\mathcal{Z}^{(j)}$. In fact, $X_j \subset [0,1]^{\mathcal{Z}^{(j)}}$. Denote $\prod_{j \in J} X_j$ by $\mathbf{X}$.

For $x = (x_1, \ldots, x_{|J|}) \in \mathbf{X}$, denote by $\mu_x^{(j)}$ the product measure on the space $\prod_{\lambda \in \Lambda_j} Z^{(\lambda)}$ such that the measure on each $Z^{(\lambda)}$ is given by $x_{\gamma(\lambda)} \in X_{\gamma(\lambda)}$.

Define an operator $\Psi_p : \mathbf{X} \to \mathbf{X}$ [where $\Psi_p(x) = (\Psi_p(x)_1, \ldots, \Psi_p(x)_{|J|})$] by

$$\Psi_p(x)_j(A) = \int_{\zeta \in \prod_{\lambda \in \Lambda_j} Z^{(\lambda)}} \mathrm{P}_p(A|\zeta) \, d\mu_x^{(j)}(\zeta) \qquad \text{for all } A \in \mathcal{Z}^{(j)}. \tag{2}$$

For any $p \in (0,1)$, the measure on $Z^{(j)}, j \in J$, induced by the percolation on $U(P_j)$ is a fixed point of this operator $\Psi_p$; indeed, compare (2) with (1).

If the pieces are finite, then the equations (1) form a finite system of polynomial equations in the unknown variables $P_p(z^{(j)})$. The probabilities of partitions are important for evaluating the first moment matrix of the branching process defined in Section 4.5.

We can endow each $Z^{(j)}$ with a topology generated by the cylindrical sets $z^{(j)}(F)$. This topological space is second countable and Hausdorff. Moreover, the space is compact since each sequence has an accumulation point. (Indeed, let $F_i, i = 1, \ldots,$ be an increasing sequence of finite subsets of $B_j$, and $\bigcup F_i = B_j$. In order to find an accumulation point of a sequence $x_k$, we can consider a set of partitions that agree on $F_i$ with infinitely many $x_k$'s and let $i$ go to infinity.)

The space **X** of probability measures on a compact metric space is a convex compact metric space in the weak-* topology induced by continuous functions (this follows from the Riesz representation theorem). Then, a continuous operator on **X** has a fixed point. The set of its fixed points is closed and since the operator $\Psi_p$ is linear, the set is also convex. The operator $\Psi_p$ acts linearly on the space of measures and has norm at most one because it preserves the subset of probability measures. Therefore, the set of its fixed points is convex, compact.

## 4. Branching processes.

4.1. *Preliminaries.* Recall that a multi-type branching process is a Markov process that models a population in which each individual in generation $n$ produces some random number of offspring of the various types in generation $n + 1$, according to a fixed probability distribution that depends only on the type of the individual.

Assume that $S$ is a set of types of individuals and that a type-$s$ individual produces children of different types according to a probability distribution $p_s$ on $\mathbb{N}_0^S$, where $\mathbb{N}_0 = \{0, 1, 2, \ldots\}$. Assume that all individuals produce offspring independently of each other and of the history of the process. Let vector $X_n \in \mathbb{N}_0^S$ represent the $n$th generation, where each coordinate $X_n[s]$ represents the number of individuals of type $s$ in the $n$th generation. It is given by the recurrent relation

$$(3) \qquad X_{n+1} = \sum_{s \in S} \sum_{i=1}^{X_n[s]} \xi_{ns}^i,$$

where $\xi_{ns}^i$ are independent random variables with the distribution $p_s$. The sequence $\{X_n\}_0^\infty$ is called a *multi-type Galton–Watson branching process* with initial population size $X_0 \in \mathbb{N}_0^S$ and offspring distribution $p_s, s \in S$.

Let $e_s \in \mathbb{N}_0^S$ be a vector with 1 at the position $s \in S$ and 0's at other positions. Let $M_n(e_s, B)$ be the expected number of individuals of the $n$th generation of types in $B \subset S$. For any fixed initial $e_s$, $M_n(e_s, \cdot)$ is a measure on $S$ given by

$$(4) \qquad M_n(e_s, B) = E\left(\sum_{s \in B} X_n[s] \Big| X_0 = e_s\right).$$

The initial population can be given by any vector $a$ in $\mathbb{N}_0^S$ (usually, it is some $e_s$, an atomic measure). We will give some properties of the $M_n$'s, using notation from the book by Nummelin [9]. Assuming that the expected size of the first generation is bounded, we can condition on the states in the first generation and obtain

$$M_n(a, B) = \int_{s \in S} M_{n-1}(s, B) \, dM_1(a, s). \tag{5}$$

Each $M_n : (S, P(S)) \to \mathbb{R}_+$ is a kernel acting as an operator $\tilde{M}_n$ on the space of measures on $S$ by

$$\tilde{M}_n(\mu)(B) = \int_{s \in S} M_n(s, B) \, d\mu(s).$$

The product of two kernels is defined as in expression (5). Thus, $\tilde{M}_n$ is the $n$th iterate kernel of $\tilde{M}_1$.

A branching process is called *singular* if each individual has exactly one offspring almost surely. $M$ is *irreducible* if and only if there exists a $\sigma$-finite measure $\phi$ on $S$ such that for all $\phi$-positive sets $B \subset S$ and $s \in S$, $M_n(e_s, B) > 0$ for some $n$.

If a nonsingular branching process is irreducible, then the population becomes extinct or explodes exponentially. The irreducibility condition is very important here. The branching process we will construct can be reducible in general. This prevents us from making claims about the extinction of the process in general. If we restrict our attention to the finite set of types, we can make further claims; see Section 4.5.

If $S$ is finite, then the operators $M_n$ are matrices. The expected size of the first generation (the first moment matrix $M$) is given by

$$m_{rs} = \mathrm{E}(X_1[s]|X_0 = e_r). \tag{6}$$

Two types $r$ and $s$ are said to be in the same class if an individual of type $r$ is in the offspring of an individual of type $s$ with positive probability and vice versa [i.e., for some $n$, the $(r,s)$ entry of $M^n$, $m_{rs}^{(n)}$, is positive]. Now, the multi-type branching process is irreducible if all types are in the same class.

The process is *positively regular* if there exists $n$ such that all elements of $M^n$ are strictly positive. If the process is irreducible but nonpositively regular, then it is said to be *periodic*. The period of a branching process is a number $d$ such that the matrix $M$ may be represented, after reordering the types of individuals, in the form

$$M = \begin{bmatrix} 0 & M(1,2) & 0 & \ldots & 0 \\ 0 & 0 & M(2,3) & \ldots & 0 \\ \vdots & \vdots & \vdots & \ddots & \vdots \\ 0 & 0 & 0 & \ldots & M(d-1,d) \\ M(d,1) & 0 & 0 & \ldots & 0 \end{bmatrix},$$

where $M(i, i+1)$ denotes a nonzero matrix.

We will need the following result from the theory of branching processes.

LEMMA 4.1. *If an irreducible multi-type Galton–Watson process with a finite number of types is nonsingular, then the population becomes extinct with probability one whenever the maximal eigenvalue of the first moment matrix $M$ is at most one.*

For the proof, see Mode [8], Theorem 7.1 on page 16 and Theorem 2.1 on page 54. Note that by the Perron–Frobenius theorem, the maximal (in absolute value) eigenvalue of $M$ will always be a nonnegative real number.

4.2. *Coloring of the tree of pieces according to the percolation.* Next, we will relate the percolation cluster size to the population of a multi-type branching process. The number of types will be as small as possible in order to simplify the computation at the cost that the population size will not match, but only approximate, the cluster size—in particular, it will be finite if and only if the percolation cluster is also finite.

Assume that there is a distinct vertex $o$ called the origin in the border set $B_0$. Each realization of a percolation gives rise to a coloring of the tree of pieces in the following way.

Consider a piece $P_i$ that is not the root. Consider a subgraph $U(P_i)^c$ of $\mathcal{G}$ induced by the edges $E(\mathcal{G}) \setminus E(U(P_i))$. Note that the intersection of $U(P_i)^c$ and $U(P_i)$ is the border set $B_i$ and that it is a cut set, by Lemma 2.3. The subgraph $U(P_i)^c$ contains the origin. Vertices of the border set $B_i$ may be connected by open paths in $U(P_i)^c$. This gives a new partition on $B_i$ and one class of the partition may be connected to the origin. The color of the piece is the data consisting of the partition and the distinct origin-connected class (which may be empty). Denote all possible colors $y^{(i)}$ by $Y^{(i)}$, that is, all colors that appear for some realization of percolation $\omega \in \Omega$. We say that a piece is white if it is of a color such that the class connected to the origin is empty. Observe that in this model, the colors of pieces depend not only on the parent, but also on the colors of its siblings and the whole subtree of their descendants. This does not give us a branching process directly, but this last difficulty is to be overcome.

Consider the maps $q'_i : \Omega \to Y^{(i)}$ assigning to each realization a color of the piece $P_i$. We can pull back the percolation probability measure $P_p$ to $Y^{(i)}$. This determines a $\sigma$-algebra $\mathcal{Y}^{(i)}$ on $Y^{(i)}$. We can identify the colors of the pieces $P_i$ and $P_{\gamma(i)}$ [$\gamma$ is the model map from part (5) of Definition 2.1] because the spaces $Y^{(i)}$ and $Y^{(\gamma(i))}$ are isomorphic.

LEMMA 4.2. *For every $i \in I$ and $v \in \Lambda_i$, $U(P_v)^c$ is covered by $P_i$, $U(P_i)^c$ and the collection of $U(P_\lambda)$ for $v \neq \lambda \in \Lambda_i$. So, the color of $P_v$ is determined* by the state of the edges in $P_i$, the color of $P_i$ and descendant partitions of border sets $B_\lambda$ for $\upsilon \neq \lambda \in \Lambda_i$.

PROOF. The argument is identical to the one used for Lemmas 2.2 and 3.1. □

Note that only a finite number of children will be nonwhite if $p < \min\{p_c(P_i)\}$ and that all children of a white piece will also be white.

4.3. *The complete branching process.* In the previous sections, we assigned to each piece $P_i$ (or to its border set $B_i$) a descending partition and a color based on the realization of the percolation process, that is, a pair $(z^{(i)}, y^{(i)}) \in \bigcup_{j \in J}(Z^{(j)} \times Y^{(j)})$. Denote $\bigcup_{j \in J}(Z^{(j)} \times Y^{(j)})$ by $\Theta$.

Assume a piece $P_j$ with border set $B_j$ has a descendant partition $z^{(j)}$ and a color $y^{(j)}$. The descendant partitions and the colors of the pieces $P_\lambda, \lambda \in \Lambda_j$ (the offspring pieces of $P_j$) are random variables taking values in $Z^{(\lambda)} \times Y^{(\lambda)}$. Denote by $(\zeta, \eta)$ the vector of these random variables, with values in $\prod_{\lambda \in \Lambda_j} Z^{(\lambda)} \times \prod_{\lambda \in \Lambda_j} Y^{(\lambda)}$. Let $Q$ be a random variable representing all the colors and descending partitions of pieces in $U(P_j)^c$.

LEMMA 4.3. *We claim that the distribution of the descendant partitions and colors of the offspring pieces depends only on the parent, in particular, that*

$$(7) \qquad P_p((\zeta, \eta) \in A | (z^{(j)}, y^{(j)}) \& Q) = P_p((\zeta, \eta) \in A | (z^{(j)}, y^{(j)}))$$

*for any measurable set $A \subset \prod_{\lambda \in \Lambda_j} Z^{(\lambda)} \times \prod_{\lambda \in \Lambda_j} Y^{(\lambda)}$.*

PROOF. From Lemmas 4.2 and 3.1, we have that $(\zeta, \eta)$ is independent of $Q$ for a given $(z^{(j)}, y^{(j)})$. Therefore, the claim follows.

In other words, the descendant partition and the color of the parent encodes every connectedness relation coming from $U(P_j)^c$. So, given the descendant partition and the color of the piece $P_i$, the descendant partitions and the colors of pieces in $U(P_j)^c$ and of pieces in $U(P_j)$ are independent. □

The purpose of this claim becomes clear after the following definition.

Let us define a measure $\mathcal{D}(z^{(i)}, y^{(i)})$ on $\prod_{\lambda \in \Lambda_i} Z^{(\lambda)} \times \prod_{\lambda \in \Lambda_i} Y^{(\lambda)}$ according to the probabilities in Lemma 4.3. We identify the space $Z^i \times Y_i$ with $Z^{\gamma(i)} \times Y^{\gamma(i)}$ and count the number of repetitions of each $(z^{(j)}, y^{(j)}), j \in J$, among the offspring. In this way, we define a map $\prod_{\lambda \in \Lambda_j} Z^{(\lambda)} \times \prod_{\lambda \in \Lambda_i} Y^{(\lambda)} \to \mathbb{N}_0^\Theta$. For any $i \in J$, we can pull back the distribution $\mathcal{D}(z^{(i)}, y^{(i)})$ to a distribution $\mathcal{D}'(z^{(i)}, y^{(i)})$ on $\mathbb{N}_0^\Theta$.

DEFINITION 4.4. The *complete* multi-type Galton–Watson branching process induced by the percolation with parameter $p$ on a graph with a tree-like structure is given by the following conditions.

The type of an individual is given by the color and the descendant partition. The set of types is $\Theta = \bigcup_{j \in J}(Z^{(j)} \times Y^{(j)})$.

There is one initial individual of type $(z_0, y_0)$, where $z_0$ is the diagonal descendant partition of $B_0$ (i.e., all classes are of size 1) and $y_0$ is the color with the diagonal partition of $B_0$ with only the origin in the distinct origin-connected class.

Every individual of some type $(z^{(j)}, y^{(j)})$ gives birth to $|\Lambda_j|$ individuals, with the distribution of types $\mathcal{D}'(z^{(j)}, y^{(j)})$.

We can represent this branching process by a tree. Clearly, it will match the tree of pieces of the graph.

LEMMA 4.5. *The coloring of the tree of pieces according to the percolation has the same distribution as the complete branching process defined above.*

PROOF. Let us represent the coloring of the tree of pieces by a random process $(X_n)$, where $X_n \in \mathbb{N}_0^\Theta$. The colors and the descendant partitions are represented by elements of $\Theta$ and each coordinate of $X_n$ gives a number of pieces with a specific color and a specific descending partition in the $n$th generation of the tree of pieces.

By Lemma 4.3, the distribution of the offspring of an individual in the $n$th generation is independent of the other individuals (in generation at most $n$) and depends only on its type $(z^{(j)}, y^{(j)})$. The distribution is $\mathcal{D}'(z^{(j)}, y^{(j)})$.

Therefore, the process $X_n$ coincides with the complete branching process above. □

4.4. *The branching process with a reduced number of types.* Next, we will reduce the number of types and obtain a different branching process. We will also impose an independence condition on the offspring of every individual, that is, the joint distribution of the offspring $\mathcal{D}(y^{(j)})$ will be a product measure. Nevertheless, we will show that the expected population size of this reduced branching process contains enough information about $p_c$.

Using Lemma 4.2, we can make the following observation about the coloring of the tree of pieces. We look at the piece $P_i$ and evaluate the conditional probability of a specific piece $P_v$ having color $y^{(v)}$ (resp., $y^{(\gamma(v))}$), assuming the descendant partitions of the other border sets are given by the vector $\zeta \in \prod_{\lambda \in \Lambda_i} Z^{(\lambda)}$ (note that the descendant partition on $P_v$ has no influence). Taking an expected value of these conditional probabilities over all possible descendant partitions gives us the probability $P_p(y^{(v)}|y^{(i)})$ of having a child piece of a given color $y^{(v)}$ from a piece of color $y^{(i)}$. These probabilities are the same for the piece $P_i$ as for its model piece $P_{\gamma(i)}$, so we can write, for $j \in J$,

$$(8) \quad P_p(y^{(\gamma(v))} \in B|y^{(j)}) = \int_{\zeta \in \prod_{\Lambda_j} Z^{(\lambda)}} P_p(y^{(\gamma(v))} \in B|y^{(j)} \& \zeta) \, d\mu^{(j)}(\zeta),$$

where $v \in \Lambda_j$ and $B \in \mathcal{Y}^{(\gamma(v))}$, a measurable set of colors. These conditional probabilities may differ for different children $v_1$ and $v_2 \in \Lambda_j$, even if $\gamma(v_1) = \gamma(v_2)$ (i.e., if the set of possible colors coincide).

Denote by $Y'^{(j)}$ the subset of $Y^{(j)}$ such that $y^{(j)} \in Y'^{(j)}$ if and only if it is not white. Let $\mathbf{Y} = \bigcup_{j \in J} Y'^{(j)}$. Then, $\mathbf{Y}$ is the set of all nonwhite colors. Let the $\sigma$-algebra $\mathcal{Y}$ on $\mathbf{Y}$ be generated by the intersections of sets in $\mathcal{Y}^{(j)}$ with $Y'^{(j)}$.

DEFINITION 4.6. The *reduced* multi-type Galton–Watson branching process induced by the percolation with parameter $p$ on a graph with a tree-like structure is given by the following conditions:

the types of individuals are colors in $\mathbf{Y}$;

there is one initial individual with color $y_0^{(0)}$, that is, the diagonal partition of $B_0$ with only the origin in the distinct origin-connected class;

every individual of some color $y^{(j)}$ gives birth to possibly $|\Lambda_j|$ individuals, each of which is born and has its color assigned independently of the others;

the distribution of the color $y^{(\gamma(\lambda))}$ of a child indexed by $\lambda \in \Lambda_j$ follows the law in formula (8).

Note that white is no longer a legitimate color. A child indexed by $\lambda \in \Lambda_j$ is not born in the reduced branching process with the same probability as that the corresponding piece is white in the percolation, that is, $P_p(y^{(\gamma(v))}$ is white $|y^{(j)})$. This is not the only difference between the reduced branching process and the coloring of the tree of pieces according to the percolation. The joint distribution of the offspring of an individual is different (because of the independence) and only the first moment is the same. However, the first moment is all that we need.

REMARK 4.7. In the definition of the tree-like structure, we allowed some isomorphism classes of pieces to be finite. In particular, very often, the root is not isometric to any other piece. In such a case, it is useful to start the branching process not at the root (with one element), but with a generation that already consists only of pieces that are in infinite isomorphism classes. This gives us the initial measure on the colors. Therefore, we consider only a subset $J'$ of $J$ such that $j \in J'$ if and only if $|\gamma^{-1}(j)| = \infty$. Then, $\mathbf{Y} = \bigcup_{j \in J'} Y'^{(j)}$.

THEOREM 1.1(i). *For a percolation with parameter $p$, the reduced branching process on the tree of pieces has the property that the expected size of its population is finite if and only if the expected size of the percolation cluster at the origin is finite.*

PROOF. Note that, by a result of Aizenman and Barsky [1], the subcritical phase on transitive graphs is equivalently characterized by the finiteness of the expected cluster size (they proved this for $\mathbb{Z}^d$, the generalization to all transitive graphs was pointed out by Lyons and Peres [7]).

Clearly, $p_c(\mathcal{G}) \leq p_c(P_i)$. If there is an infinite cluster at the origin for $p$ such that $p < p_c(P_i)$, then the cluster must intersect infinitely many border sets almost surely. Therefore, the number of border sets connected to the origin is finite if and only if the original cluster was also finite. In what follows, we will always assume that $p \leq p_c(P_i)$ [if $P_i$ is finite, then we set $p_c(P_i) = 1$].

Given a realization of percolation, we have introduced a coloring of the tree of pieces. The percolation cluster at the origin is infinite if and only if the nonwhite colored component of the tree of pieces is infinite (by the assumption that $p$ is smaller than $p_c$ of the pieces). The expected size of the colored component is a sum of the probabilities that a piece is nonwhite over all pieces of the tree.

The branching process from Definition 4.6 can be naturally illustrated by a tree isomorphic to the tree of pieces. The distribution of colors of a specific individual in this branching process equals the distribution of nonwhite colors of the related piece in the percolation. Therefore, the expected population size of this branching process equals the expected number of nonwhite colored pieces in the percolation.

Therefore, $p < p_c$ if and only if the related branching process has finite expected population size. □

The branching process is not singular because, for $p < 1$ (and some nonwhite color $y$), there is no color which appears as an only child of $y$ almost surely (i.e., with positive probability, there are more or less children with different colors). Assume that the initial measure is $\nu : \mathcal{Y} \to \mathbb{R}_+$ (it can be the atomic measure, from Definition 4.6, or the distribution of the first generation of pieces with infinite isomorphism classes, from Remark 4.7). Then, the expected population size of the branching process is $\sum_{n=0}^{\infty} \tilde{M}_n \nu(\mathbf{Y})$. The critical probability $p_c$ is then a supremum of all $p$ such that $\sum_{n=0}^{\infty} \tilde{M}_n \nu(\mathbf{Y}) < \infty$. In general, the convergence of $\sum_{n=0}^{\infty} \tilde{M}_n \nu(\mathbf{Y})$ may depend strongly on the initial measure $\nu$. However, we will show that if the border sets are finite, then it actually depends only on the spectral radius of $M_1$.

4.5. *The case of finite border sets.* If the border sets are finite, then the space of partitions is finite and so is the space of colors **Y**. Denote by $M = [m_{ab}]_{a,b \in \mathbf{Y}}$ the first moment matrix of the branching process, that is, the matrix of expected numbers of offspring of each color $m_{ab} = \mathrm{E}_p(\#b|a)$, where

$$\mathrm{P}_p(y^{(\gamma(v))}|y^{(j)}) = \int_{\zeta \in \prod_{\Lambda_j} Z^{(\lambda)}} \mathrm{P}_p(y^{(\gamma(v))}|y^{(j)} \& \zeta) \, d\mu^{(j)}(\zeta),$$

(9)

$$\mathrm{E}_p(\#y^{(k)}|y^{(j)}) = \sum_{\lambda \in \Lambda_j, \gamma(\lambda) = k} \mathrm{P}_p(y^{(\gamma(\lambda))}|y^{(j)}).$$

The expected number of individuals of the $n$th generation is then given by the $n$th power of $M$.

THEOREM 1.1(ii). *If all of the border sets are finite, then the branching process has finitely many types and the first moment matrix is of finite size. In this case, $p_c$ is the smallest value of $p$ such that the spectral radius of the first moment matrix is 1.*

PROOF. Assume that the graph $\mathcal{G}$ has a tree-like structure such that the border sets $B_i$ are finite. Using the already-proven part (i) of Theorem 1.1, we need to decide for which $p$ the expected population size of the constructed branching process is finite.

This branching process is nonsingular for $p < 1$ because if a piece has more than one child, then the offspring size is greater than one with positive probability, and if every piece has exactly one child, then it has no offspring (only white) with positive probability. The expected population size is $\sum M^n$ applied to the initial measure. If the spectral radius is less than 1, then the sum $\sum M^n$ is always finite. If the spectral radius is at least 1, then there exists a possible initial measure, for which the expected population size is infinite (this follows from the Perron–Frobenius theorem). If the process is irreducible, then it is independent of the choice of initial measure.

In the case where $M$ is irreducible, the result follows directly from Lemma 4.1. In the other case, there are several classes of types and we denote by $M_1, \ldots, M_k$ the first moment matrices of each class. We can reorder the types in **Y** such that $M$ has diagonal blocks equal to the $M_i$'s and all entries above these diagonal blocks are zero. Thus, the spectral radius of $M$ is the spectral radius of some $M_s$. There is a type (color from **Y**) in the $s$th class such that, with positive probability, there is a piece $P_i$ with this color (because we considered only those colors that are realized by percolation). The process starting at this $P_i$ will have infinite cluster size whenever the spectral radius of $M_s$ is at least one. Therefore, the expected size of the percolation cluster is finite whenever the maximal eigenvalue of $M$ is less than one.

Note that $\det(M - 1) = 0$ if 1 is an eigenvalue of $M$, and for $p = 0$, all eigenvalues of $M$ are zero. The eigenvalues depend continuously on the matrix entries, which are continuous functions of $p$. Therefore, $p_c$ is the first positive value of $p$ such that $\det(M - 1) = 0$. □

THEOREM 1.1(iii). *If all of the pieces are finite, then the entries of the first moment matrix of the reduced branching process are algebraic functions in $p$. Therefore, $p_c$ is algebraic.*

*There exists an algorithm that, given the model pieces and their border sets, computes a finite extension $K$ of the field $\mathbb{Q}(p)$ and a function $f$ in $K$ such that $p_c$ is the smallest positive root of $f$.*

PROOF. If the pieces are finite, then the probabilities of the descending partitions, $P_p(z^{(j)})$, can be found as solutions of the system of equations (1) introduced in the previous section. Let us index the partitions on $B_j$ by natural numbers $1, \ldots, |Z^{(j)}|$. Denote by $x_{j,i}$ the probability of the $i$th possible partition on the border set of model piece $P_j$, that is, $x_{j,i} = P_p(z_i^{(j)})$. We then have the following system of equations in the unknown variables $x_{j,i}$:

$$\tag{10} x_{j,i} = \sum_k \left( \sum_{\Gamma \in L} \alpha(\Gamma) \right) \prod_{\lambda \in \Lambda_j} x_{\gamma(\lambda), k(\lambda)},$$

where the following statements hold true:

- $k$ is a map that assigns partitions to the border sets $B_\lambda$. That is, $k: \Lambda_j \to \mathbb{N}$, and for all $\lambda \in \Lambda_j$, $k(\lambda) \in \{1, 2, \ldots, |Z^{(\lambda)}|\}$. We sum over all such possible functions $k$.
- For a given $k$, $L$ is the set of subgraphs of $P_j$ such that $\Gamma \in L$ if and only if the following holds: provided the partitions on the children pieces are $z_{k(\lambda)}^{(\lambda)}$ and the open edges of $P_j$ are given by $\Gamma$, the partition on $B_j$ [induced by $U(P_j)$] equals $z_i^{(j)}$. Note that, by Lemma 3.1, the partition on $B_j$ is uniquely determined by the information provided.
- $\alpha(\Gamma)$ is the probability of $\Gamma$ in the percolation, that is, $\alpha(\Gamma) = p^{E(\Gamma)}(1 - p)^{E(P_i) - E(\Gamma)}$.

The number of equations is equal to the number of variables. The degree in $p$ of each equation (10) is equal to the number of edges in a piece and the degree in the unknown variables $x_{j,i}$ is $|\Lambda_j|$, the number of offspring of the piece $P_j$. Denote by $K$ the algebraic extension of $\mathbb{Q}(p)$ containing roots of this system of equations.

Similarly as above, we can rewrite the formulae (9) as follows:

$$P(y^{(\gamma(v))}|y^{(j)}) = \sum_{k}\left(\sum_{\Gamma \in L'} \alpha(\Gamma)\right) \prod_{\lambda \in \Lambda_j} x_{\gamma(\lambda),k(\lambda)}, \tag{11}$$

where:

- $k$ and $\alpha(\Gamma)$ are defined as above, and
- $L'$ depends on $k$ and it is the set of subgraphs of $P_j$ such that $\Gamma \in L'$ if and only if the color on $B_v$ is $y^{(\gamma(v))}$ provided the color of the parent is $y^{(j)}$, the descendant partitions on the children pieces are $z_{k(\lambda)}^{(\lambda)}$ and the open edges of $P_j$ are given by $\Gamma$ (by Lemma 4.2, these uniquely determine the color on $B_v$).

Consequently, the entries of the first moment matrix $M$, given by formula (9), are polynomial functions in $p$ and the $x_{j,i}$'s. The critical probability $p_c$ is the first positive value of $p$ such that $\det(M-1) = 0$. The function $\det(M-1)$ is an element of the field $K$ and $p_c$ is an algebraic number in this case. □

## 5. Examples.

5.1. *Free products of (transitive) graphs.* The free product $\mathcal{G}$ of transitive graphs $G_1, \ldots, G_n$ is an infinite connected graph constructed as a union of copies of $G_i$'s such that each vertex of $\mathcal{G}$ belongs to exactly one copy of each $G_i$ and every simple closed path in $\mathcal{G}$ is included in one copy of some $G_i$. In particular, if $G_1 = \langle S_1 \rangle$ and $G_2 = \langle S_2 \rangle$, then a Cayley graph of $G_1 * G_2$ with respect to $S_1 \cup S_2$ is a free product of Cayley graphs of $G_1$ and $G_2$ with respect to $S_1$ and $S_2$, respectively.

Here, we present a different method to obtain a result shown previously by the author in [6].

COROLLARY (1.2). *Let $\mathcal{G}$ be a free product of transitive graphs $\mathcal{G}_1, \ldots, \mathcal{G}_n$. Denote by $\chi_i(p)$ the expected subcritical cluster size in the ith factor graph $G_i$. The critical probability $p_c$ of $\mathcal{G}$ is the infimum of positive solutions of $\sum_{j=1}^{n} \prod_{i=1, i\neq j}^{n} \chi_i(p) - (n-1) \prod_{i=1}^{n} \chi_i(p) = 0$.*

PROOF. The tree-like structure is very natural in this case. The copies of factor graphs $\mathcal{G}_1, \ldots, \mathcal{G}_n$ are pieces (thus, $J = \{1, \ldots, n\}$) and each vertex is, in fact, a border set.

The partitions of border sets are trivial (the spaces of partitions have size 1) in this case. For the coloring, we need only two colors: $Y^{(i)} = \{y_1^{(i)}, y_2^{(i)}\}$. Assuming that $y_2^{(i)}$ is white, we have $\mathbf{Y} = \{y_1^{(1)}, \ldots, y_1^{(n)}\}$.

Denote by $b_i$ the border vertex of piece $P_i$ and by $\tau_p^{(i)}(a \leftrightarrow b)$ the probability that $a$ and $b$ are connected in $\mathcal{G}_i$. Denote by $\chi_i(p)$ the expected cluster size in $\mathcal{G}_i$. We are now ready to compute the matrix $M = (m_{ij})$:

$$P_p(y_1^{(\gamma(\lambda))}|y_1^{(j)}) = \tau_p^{(j)}(b_j \leftrightarrow b_\lambda),$$
$$E_p(\#y_1^{(k)}|y_1^{(j)}) = \chi_j(p) - 1,$$
$$m_{ij} = \chi_i(p) - 1 \quad \text{for } i \neq j,$$
$$m_{ii} = 0.$$

In order to find for which $p$ the spectral radius of $M$ equals 1, we solve the equation $\det(M - 1) = 0$ by Theorem 1.1(iii). The determinant is computed as follows (we subtract the first column from all of the others and then expand the determinant along the first row):

$$\det(M - 1) = \det \begin{bmatrix} -1 & \chi_1(p)-1 & \cdots & \chi_1(p)-1 \\ \chi_2(p)-1 & -1 & \cdots & \chi_2(p)-1 \\ \vdots & \vdots & \ddots & \vdots \\ \chi_n(p)-1 & \chi_n(p)-1 & \cdots & -1 \end{bmatrix}$$

$$= \det \begin{bmatrix} -1 & \chi_1(p) & \chi_1(p) & \cdots & \chi_1(p) \\ \chi_2(p)-1 & -\chi_2(p) & 0 & \cdots & 0 \\ \chi_3(p)-1 & 0 & -\chi_3(p) & \cdots & 0 \\ \vdots & \vdots & \vdots & \ddots & \vdots \\ \chi_n(p)-1 & 0 & 0 & \cdots & -\chi_n(p) \end{bmatrix}$$

(12)
$$= -\prod_{i=2}^{n}(-\chi_i(p)) + \sum_{j=2}^{n}(\chi_j(p)-1)\prod_{i=1, i\neq j}^{n}(-\chi_i(p))$$

$$= (-1)^n \left( \sum_{j=1}^{n} \prod_{i=1, i\neq j}^{n} \chi_i(p) - (n-1)\prod_{i=1}^{n} \chi_i(p) \right).$$

Therefore, $p_c$ is the infimum of positive solutions of

$$\sum_{j=1}^{n} \prod_{i=1, i\neq j}^{n} \chi_i(p) - (n-1)\prod_{i=1}^{n} \chi_i(p) = 0.$$

$\square$

5.2. $SL(2,\mathbb{Z})$. The easiest example of a Cayley graph with a tree-like structure that is not a free product is $SL(2,\mathbb{Z})$. It is the amalgamated product $\mathbb{Z}_4 *_{\mathbb{Z}_2} \mathbb{Z}_6$ with standard generating set $\{a,b\}$ so that $\mathbb{Z}_4 = \langle a \rangle$ and $\mathbb{Z}_6 = \langle b \rangle$. We can illustrate the general method here because the required computations (for obtaining $p_c$) are relatively simple, or at least doable.

COROLLARY (1.3). *The critical probability $p_c$ of the special linear group $SL(2, \mathbb{Z})$ given by presentation $\langle a, b | a^4, b^6, a^2 b^{-3} \rangle$ is an algebraic number equal to $0.4291140496\ldots$.*

PROOF. The Cayley graph has the tree-like structure with pieces of two isomorphism classes (squares $P_1$ and hexagons $P_2$) corresponding to the factor groups; see Figure 1. The hexagonal pieces consist of six-tuples of vertices connected by bold lines in Figure 1 and squares contain dashed edges. Each vertex is contained in exactly one square and one hexagon, and each neighboring square and hexagon share two vertices—the border sets.

First, we need to find the distribution of the descendant partitions of the border sets. There are two possible descendant partitions on each border set, $B_1$ and $B_2$. Denote them by $Z^{(1)} = \{z_1^{(1)}, z_2^{(1)}\}$ and $Z^{(2)} = \{z_1^{(2)}, z_2^{(2)}\}$, where $z_1^{(i)}$ means that the two vertices of the border set $B_i$ are connected and $z_2^{(i)}$ means that they are not connected. We can now find the probability of these partitions using formula (1), which leads to a quadratic equation.

In Section 4, we introduced the coloring of a tree. For a border set of size two, there are three nonwhite colors: either both vertices are connected to the origin or exactly one of them is connected to the origin. The two different situations when exactly one vertex of the border set is connected to the origin are symmetric, that is, there is an obvious isomorphism of $U(P_i)$ to itself such that each vertex of $B_i$ is mapped onto the other one. Therefore, these two colors have the same distribution of offspring. This simplification is special for the Cayley graph of $SL(2, \mathbb{Z})$ because this new symmetry is not one of the isomorphisms in condition (5) of the definition of the tree-like structure. Since the distribution of offspring does not depend on which of the two vertices is connected to the origin, we identify these two colors as one.

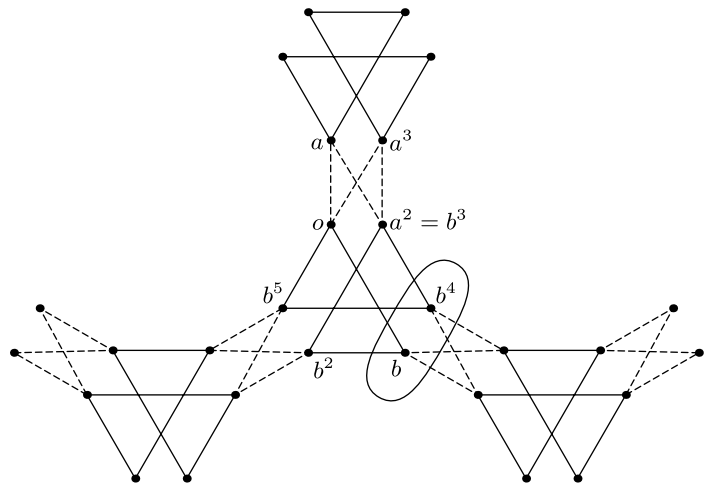

FIG. 1. *The Cayley graph of $SL(2, \mathbb{Z})$ with one border set circled.*

The set of colors becomes $Y^{(i)} = \{y_r^{(i)}, y_b^{(i)}, y_w^{(i)}\}$, where $r$ (red) means that both vertices of the border set are connected to the origin, $b$ (blue) means that exactly one vertex is connected to the origin (thus the border set is disconnected) and $w$ stands for white—no vertex connected to the origin. We evaluate the probabilities of getting a child of a given color, according to formula (11), and the first moment matrix $M$ takes the following form:

$$(13) \quad M = \begin{bmatrix} 0 & 0 & P_p(y_r^{(\lambda)}|y_r^{(1)}) & P_p(y_b^{(\lambda)}|y_r^{(1)}) \\ 0 & 0 & P_p(y_r^{(\lambda)}|y_b^{(1)}) & P_p(y_b^{(\lambda)}|y_b^{(1)}) \\ 2P_p(y_r^{(\lambda)}|y_r^{(2)}) & 2P_p(y_b^{(\lambda)}|y_r^{(2)}) & 0 & 0 \\ 2P_p(y_r^{(\lambda)}|y_b^{(2)}) & 2P_p(y_b^{(\lambda)}|y_b^{(2)}) & 0 & 0 \end{bmatrix}.$$

If we solve the equation $\det(M-1) = 0$, we obtain $p_c = 0.4291140496\ldots$.
□

5.3. *Grandparent tree.* An interesting example of a transitive graph with infinitely many ends that is not a Cayley graph is the grandparent tree. It is obtained from the three-regular tree, by adding edges connecting each vertex to its grandparent, that is, we pick one end and, at each vertex, we add an edge connecting it to a vertex at distance 2 in the direction of the distinguished end.

In the suggested decomposition below, we obtain pieces of size four with border sets of size two. Note that in Figure 2, the graph is oriented so that the distinguished end is at the top. In the right-hand figure, we see that the border set of a piece consists of the two upper vertices and each of its two children intersects it by a pair of vertices—the middle one and the left (resp., right) one.

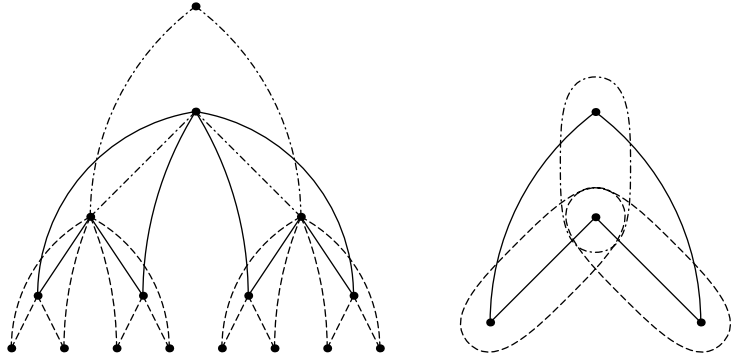

FIG. 2. *Part of the grandparent tree and its piece with cycled border sets.*

Note that the pieces that lie on the path from the root to the distinguished end have a different orientation—the border set is a different pair of vertices inside the piece (left as an exercise for the reader).

A careful analysis of the branching process can be carried out as in the case of $SL(2,\mathbb{Z})$. The $p_c$ obtained is $0.158656326\ldots$.

**6. Fundamental groups of graphs of groups.** In this section, we will generalize the example of $SL(2,\mathbb{Z})$ to arbitrary graphs of groups by showing that the Cayley graphs of fundamental groups of graphs of groups have the tree-like structure. To make the transition simpler, we first consider amalgamated products and HNN extensions before proceeding to arbitrary graphs of groups.

6.1. *Amalgamated products with standard generating sets.* Recall that an amalgamated product $G_1 *_H G_2$ is obtained from two groups $G_1$ and $G_2$ provided with monomorphisms $i_k : H \to G_k$, $k = 1, 2$. The group $G_1 *_H G_2$ is a quotient of the free product $G_1 * G_2/N$, where $N$ is the smallest normal subgroup containing elements $i_1(h)i_2^{-1}(h)$, $h \in H$.

Consider the right Cayley graph $\mathcal{G}$ of an amalgamated product $G_1 *_H G_2$ with respect to the generating set $S_1 \cup S_2$, where $G_1 = \langle S_1 \rangle$ and $G_2 = \langle S_2 \rangle$. The structure of the Cayley graph $\mathcal{G}$ is the following. First, consider the copies of the Cayley graph $\mathcal{G}_1$ of $G_1$ and the Cayley graph $\mathcal{G}_2$ of $G_2$ containing the origin. These subgraphs intersect at the vertices of $H$, which we consider as a border set. Each coset of $H$ inside $G_1$ is another border set connecting $\mathcal{G}_1$ with another copy of $\mathcal{G}_2$; each coset of $H$ inside $G_2$ is a border set connecting it with another copy of $\mathcal{G}_1$. The copies of $\mathcal{G}_1$ and $\mathcal{G}_2$ form a tree—it is a $G$-tree (see [2], Section 8.6) which is usually called the *Bass–Serre tree* of the amalgamated product.

Let us pick one of the pieces, say $\mathcal{G}_1$, and call it a root $P_0$. It contains $[G_1 : H]$ border sets $B_\lambda$, $\lambda \in \Lambda_1$. All of these $P_\lambda$, $\lambda \in \Lambda_1$, are in the same isomorphism class, say $I_2$, because there exists an isomorphism that takes $P_\lambda$ to $P_{\lambda'}$ with $\lambda, \lambda' \in \Lambda_1$ and carries descendants to descendants. This isomorphism is just a left multiplication by an element of $G_1 *_H G_2$. Children of pieces in the isomorphism class $I_2$ are in the same isomorphism class again—denote it by $I_1$. So, the isomorphism class of a piece depends only on the parity of the generation. More precisely, there are three isomorphism classes: the root, the set of pieces of odd generations and the set of pieces of even generations. The root differs from an even generation piece by the number of children, but it is only one such exceptional piece. No color assigned to this piece is in $\mathbf{Y}$, by Remark 4.7, so we do not need to consider this piece in the remaining computations.

The parity of a generation implies that the branching process defined in Section 4 is periodic and so the matrix $M$ always has two anti-diagonal blocks, as in the case of $SL(2,\mathbb{Z})$; see (13).

6.2. *HNN extensions with standard generating sets.* Recall that an HNN extension $G$ is constructed from a base group $G_1$ having a presentation $G_1 = \langle S|R \rangle$ and from an isomorphism $\alpha$ between two subgroups $H$ and $K$ of $G_1$. Let $t$ be a new symbol not in $S$ (free letter) and define $G = G_1 *_\alpha = \langle S, t | R, tht^{-1} = \alpha(h), \forall h \in H \rangle$.

Consider the Cayley graph $\mathcal{G}$ of this HNN extension $G$ with respect to the generating set $S \cup \{t\}$. The Cayley graph $\mathcal{G}_1$ of $G_1$ is a part of the Cayley graph $\mathcal{G}$. Each coset of $H$ (resp., $K$) is attached to another copy of $\mathcal{G}_1$—the attachment is done by edges labeled by $t$, these edges corresponding to the isomorphism $\alpha$. Therefore, the piece $\mathcal{G}_1$ is connected to $[G_1:H] + [G_1:K]$ other pieces. In this way, we obtain the Bass–Serre tree of the HNN extension.

Let us denote by $P_0$ a subgraph of $\mathcal{G}$ containing the graph $\mathcal{G}_1$ and all $t$-edges incident to at least one vertex in $\mathcal{G}_1$ (and we also add its other endpoint). In $\mathcal{G}$, the origin is connected by a $t$-edge to some vertex $v$, we denote by $P_1$ a subgraph of $\mathcal{G}$ that contains the copy of $\mathcal{G}_1$ at $v$ and all $t$-edges incident to at least one vertex in this copy that are not in $P_0$. The border set $B_1$ consists of vertices in $H$ in the copy of $\mathcal{G}_1$. Repeating this for all $t$-edges emerging from $P_0$, we obtain $[G_1:K]$ pieces in the first generation. Similarly, by following edges labeled by $t^{-1}$ emerging from $P_0$, we obtain $[G_1:H]$ other pieces in the first generation. We then do the same for the next generations. This procedure gives the tree-like structure of the graph. Each piece $P_i$ contains a copy of $\mathcal{G}_1$ and some $t$-edges connected to it (in fact, it contains all such $t$-edges except the $|K|$ edges that connect the border set $B_i$ and the parent). There are three isomorphism classes of pieces: the root, pieces whose border sets are copies of $K$ and pieces whose border sets are copies of $H$. In this case, the process is not periodic since a piece of each class may be incident to some other piece of the same class.

6.3. *Graphs of groups.* Let us recall the definition of the fundamental group of a graph of groups [2].

DEFINITION 6.1. A *graph of groups* $\mathfrak{G}$ consists of:

(i) a connected graph $X$ with vertex set $V(X)$ and edge set $E(X)$;
(ii) for each vertex $v$ of $X$, a group $G_v$, and for each edge $e$ of $X$, a group $G_e$;
(iii) for each edge $e = (v_1, v_2)$, monomorphisms $\tau : G_e \to G_{v_1}$ and $\sigma : G_e \to G_{v_2}$.

Denote by $E$ the free group with basis $\{t_e; e \in E(X)\}$. Let $F(\mathfrak{G})$ be the group $(E * *_{v \in V(X)} G_v)/N$, where $N$ is the normal closure of the subset $\{t_e^{-1} \tau(a) \times t_e \sigma(a)^{-1} : e \in E(X), a \in G_e\}$.

Let $T$ be a maximal tree in $X$. We define the *fundamental group* $\pi(\mathfrak{G}, X, T)$ to be $F(\mathfrak{G})/M$, where $M$ is the normal closure of $\{t_e, e \in E(T)\}$.

Note that the groups $\pi(\mathfrak{G}, X, T)$ are independent of $T$ up to isomorphism.

We will consider the Cayley graph $\mathcal{G}$ of the fundamental group $\pi(\mathfrak{G}, X, T)$ with the following generators: $(\bigcup_{v \in V(X)} S_v) \cup \{t_e : e \in E(X) \setminus E(T)\}$, where $G_v = \langle S_v \rangle$ for $v \in V(X)$. This set of generators depends on $T$. We call this set of generators *standard*. In order to obtain a locally finite graph, we will consider only finite sets of generators. In particular, this restricts us to finite graphs $X$ of finitely generated groups.

Observe that if $X$ consists of one edge between two distinct vertices, then its fundamental group is an amalgamated product of the vertex groups. If the two vertices coincide, then the fundamental group is an HNN extension. The tree-like structure we described for these specific cases will now be generalized to the Cayley graph of any (finite) graph of groups.

THEOREM (1.4). *The Cayley graph of the fundamental group of a graph of groups has a tree-like structure.*

PROOF. We will define the pieces and verify the conditions (1)–(5) in Definition 2.1.

First, observe that the Cayley graph $\mathcal{G}$ is covered by translates of the Cayley graphs $\mathcal{G}_v$ of the vertex groups $G_v$, $v \in V(X)$, by left multiplications by elements of $G$ and edges labeled by $t_e$, $e \in E(X) \setminus E(T)$. For each $e = (v_1, v_2) \in E(X) \setminus E(T)$, denote by $R_e$ the set of edges labeled by $t_e$ such that they start at a vertex of $\mathcal{G}_{v_1}$ and end at a vertex of the translate of $\mathcal{G}_{v_2}$ by $t_e$. We then denote by $g\mathcal{G}_v$ the translate of the subgraph $\mathcal{G}_v$ by $g$, the representative of a coset in $\pi(\mathfrak{G}, X, T)/G_v$. Denote by $ghR_e$ the translate of the set $R_e$ by $gh$, where $g$ is a representative of a coset in $\pi(\mathfrak{G}, X, T)/G_{v_1}$ and $h$ is a representative of a coset in $G_{v_1}/\tau(G_e)$. These sets will be used in the construction of pieces of the tree-like structure.

Consider the quotient map $\pi(\mathfrak{G}, X, T) \to \pi_1(X)$ taking all $G_v$ to the identity. The fundamental group $\pi_1(X)$ of a finite graph $X$ is free of rank $|E(X) \setminus E(T)|$. Every translate of $R_e$ is a preimage of one edge in the Cayley graph of $\pi_1(X)$, thus, it is an edge cut set of $\mathcal{G}$. Moreover, if $e = (v_1, v_2) \in E(T)$, then any translate $g\tau(G_e)$ by $g \in \pi(\mathfrak{G}, X, T) \to \pi_1(X)$ is a vertex cut set.

We say that $g_1\mathcal{G}_{v_1}$ and $g_2\mathcal{G}_{v_2}$ are *neighbors* if there is an edge $e = (v_1, v_2) \in E(T)$ [or $e = (v_2, v_1)$] and $g_1\mathcal{G}_{v_1} \cap g_2\mathcal{G}_{v_2} \neq \varnothing$. We say that $g_1\mathcal{G}_{v_1}$ and $g_2\mathcal{G}_{v_2}$ are *delayed neighbors* if there is an edge $e = (v_1, v_2) \in E(X) \setminus E(T)$ [resp., $e = (v_2, v_1)$], and if there is translate $hR_e$ of a set of edges labeled by $t_e$ connecting a vertex of $g_1\mathcal{G}_{v_1}$ to a vertex of $g_2\mathcal{G}_{v_2}$, then we say the neighbors are *delayed by $hR_e$*.

Consider a graph with vertices $g\mathcal{G}_v$ and edges between neighbors and delayed neighbors. We claim that this graph is a tree. Assume for a contradiction that it is not a tree and let there be a cycle. If some edge of the cycle arises from delayed neighbors, then this $R_e$ is not a cut set. If all edges in the cycle are neighbors, then they induce a loop in $T$, which is also a contradiction. We choose a root of this tree to be $\mathcal{G}_{v_o}$ for some $v_0 \in V(X)$.

We are now ready to define pieces and border sets. We take the tree from above and define the pieces: for the root $\mathcal{G}_{v_0}$, we define a piece $P_{v_0,o}$ to be a union of $\mathcal{G}_{v_0}$ and sets $hR_e$ for all its delayed neighbors. For any other $g\mathcal{G}_v$, we define a piece $P_{v,g}$ to be a union of $\mathcal{G}_{v_0}$ and sets $hR_e$ for all its delayed neighbors that are not parents of $g\mathcal{G}_v$ in the above tree order. Define the border sets of a piece $P_{v,g}$ to be its intersection with the parent piece. For the root, define the border set to contain only the origin.

Let us now verify the conditions in Definition 2.1:

(1) The Cayley graphs $\mathcal{G}_v$ are connected and the pieces remain connected after adding the incident edges $t_e$. The intersection of $P_{v_1,g_1}$ with its parent $P_{v_2,g_2}$ is a translate of $\tau(G_e)$ [resp., $\sigma(G_e)$] inside $P_{v_1,g_1}$, provided $e = (v_1, v_2)$ [resp., $e = (v_2, v_1)$]. Thus, the border set is indeed a subset of the vertex set of the piece.
(2) We need to show that every edge $e \in E(\mathcal{G})$ is in exactly one $P_i$. Clearly, this holds for edges labeled by $s \in S_v$. Every edge labeled by $t_e$ has endpoints in two different pieces which, by the construction, become parent and child, and thus it is only included in the parent piece.
(3) The tree was defined above in such a way that the border sets satisfy required condition (3).
(4) For all $j \neq 0$, the border set $B_j$ is a cut set. Using the notation $U(P_i)^c$ for the subgraph induced on edges in $E(G) \setminus E(U(P_i))$, we obtain $U(P_i) \cap U(P_i)^c = B_i$. So, claim (4) follows.
(5) We now want to show that there are only finitely many isomorphism classes of pieces. The root is a unique piece of its class and we exclude it from the following consideration. Assume that $P_i$ was obtained from a translate of $\mathcal{G}_{v_1}$ and its parent from $\mathcal{G}_{v_2}$ that are neighbors due to an edge $e = (v_1, v_2) \in E(X)$. This edge, together with its direction, characterizes the isomorphism class of $P_i$. Indeed, if $P_j$ is characterized by the same edge $e = (v_1, v_2) \in E(X)$, then, clearly, there is an element $g \in \pi(\mathfrak{G}, X, T)$ acting on the Cayley graph as an isomorphism $f$ such that $f(P_i) = P_j$ and $f(B_i) = B_j$. Therefore, $f(U(P_i)) = U(P_j)$ as well. Moreover, since the pieces (and border sets) are identified in each subtree by the same procedure, this isomorphism respects the structure of pieces in the subtree. □

COROLLARY 6.2. *For the Cayley graph of a fundamental group of a finite graph of finite groups, with the natural generating sets, the $p_c$ is an algebraic number.*

PROOF. This follows from Theorems 1.1(iii) and 1.4. □

REMARK 6.3. A fundamental group of a graph of groups with trivial edge groups is a free product of its vertex groups and copies of $\mathbb{Z}$, so the first moment matrix $M$ can be expressed as in Example 5.1.

6.4. *An example with infinite pieces.* Consider an amalgamated product $G_1 *_{\mathbb{Z}_2} G_2$ with $G_1 = \mathbb{Z}_2 \times \mathbb{Z}$ and $G_2 = \mathbb{Z}_4$. Its Cayley graph with respect to natural generators can be decomposed similarly as above, with border sets of size two. Therefore, the set of partition $Z^{(i)}$ has two elements and the set of colors **Y** has four elements, as for the Cayley graph of $SL(2, \mathbb{Z})$. The first moment matrix $M$ again has two anti-diagonal blocks.

Let $A$ (resp., $B$) be the probability that the two vertices of the border set in $G_1$ (resp., $G_2$) are connected (in the descendant subtree). $A$ and $B$ then satisfy the following:

$$A = 2p - p^2 + (1-p)^2(2C - C^2),$$
$$B = (2p - p^2)^2 A + (2p^2 - p^4)(1 - A),$$
$$C = p^2(1 - (1-p)^2(1-B) + (1-p)^2(1-B)C),$$
$$C = \frac{p^2(1 - (1-p)^2(1-B))}{1 - p^2(1-p)^2(1-B)}.$$

In order to express $A$, we split the piece $\mathbb{Z}_2 \times \mathbb{Z}$ into three parts, the first containing only two edges between vertices of $B_i$. We then split the rest at the $B_i$ and obtain two other (identical) parts (corresponding to $\mathbb{Z}_2 \times \mathbb{N}$; see, e.g., the last block in Figure 3). The probability that vertices of $B_i$ are connected in one of these two latter parts is equal to $C$.

These formulae lead to a cubic equation and exactly one of the solutions is real and in $[0, 1]$. We are now ready to express the entries of the first moment matrix $M$. There are two anti-diagonal blocks $M_{12}$ and $M_{21}$; again, the one corresponding to $\mathbb{Z}_4$ is already known from a previous example. For the other factor group $\mathbb{Z}_2 \times \mathbb{Z}$, we denote by $T_n$ the transition matrix from the parent to the $n$th child (order them by distance from the origin). In order to obtain the transition matrix, we split the graph of $\mathbb{Z}_2 \times \mathbb{Z}$ into several blocks; see Figure 3. For each block, we express the transition probabilities and thus $T_n$ is the composition of these probabilities; in particular, it is a product of matrices corresponding to each block.

The first block consists of vertices to the left of the parent and including the parent. The probability of the connection in the left part is $C$ and the probability that the two vertices in $B_0$ are connected is $2p - p^2$. This gives the first matrix in formula (14). The second block is used to express matrix $T$, that is, the transition to the first child $B_1$ from the parent $B_0$. The third one is repeated $n-1$ times in order to reach the $n$th child. It differs from $T$ by the factor of $B$ coming from the descendant partition on the sibling border set (this connection is drawn by the dashed/dotted line). The last factor in the product corresponds to the rightmost part of the picture and equals the probability of connection on the right side of the child in question (again, it is equal to $C$):

$$
\begin{aligned}
T_n &= \begin{bmatrix} \mathrm{P}_p(r_n|r) & \mathrm{P}_p(b_n|r) \\ \mathrm{P}_p(r_n|b) & \mathrm{P}_p(b_n|b) \end{bmatrix} \\
&= \begin{bmatrix} 1 & 0 \\ 1-(1-C)(1-p)^2 & (1-C)(1-p)^2 \end{bmatrix} \\
&\quad \times T\left(\begin{bmatrix} 1 & 0 \\ B & 1-B \end{bmatrix} T\right)^{n-1} \begin{bmatrix} 1 & 0 \\ C & (1-C) \end{bmatrix},
\end{aligned}
$$
(14)

$$
\begin{aligned}
T &= \begin{bmatrix} \mathrm{P}_p(r|r) & \mathrm{P}_p(b|r) \\ \mathrm{P}_p(r|b) & \mathrm{P}_p(b|b) \end{bmatrix} \\
&= \begin{bmatrix} p^2 + 2p(1-p)(2p-p^2) & 2p(1-p)^3 \\ p(2p-p^2) & p(1-p)^2 \end{bmatrix},
\end{aligned}
$$

$$
M_{21} = \begin{bmatrix} \mathrm{E}_p(r|r) & \mathrm{E}_p(b|r) \\ \mathrm{E}_p(r|b) & \mathrm{E}_p(b|b) \end{bmatrix} = 2 \sum_{n=1}^{\infty} T_n
$$

$$
= 2 \begin{bmatrix} 1 & 0 \\ C & (1-C) \end{bmatrix} \begin{bmatrix} 1 & 0 \\ 2p-p^2 & (1-p)^2 \end{bmatrix}
$$
(15)
$$
\times T\left(1 - \begin{bmatrix} 1 & 0 \\ B & 1-B \end{bmatrix} T\right)^{-1} \begin{bmatrix} 1 & 0 \\ C & (1-C) \end{bmatrix}.
$$

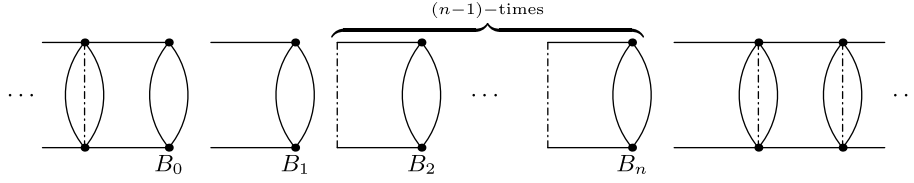

FIG. 3. *Blocks of the Cayley graph of $\mathbb{Z}_2 \times \mathbb{Z}$.*

Again, we can solve the equation $\det(M-1) = 0$ and we obtain $p_c = 0.2951\ldots$.

6.5. *Groups acting on trees with finite vertex stabilizers.* The goal of this section is to prove Theorem 1.5. To this end, we combine the result about standard generating sets (Theorem 1.4) with the general results from Section 2.

THEOREM (1.5). *Let $G$ be a virtually free group, that is, it acts on a simplicial tree $T$ with finite vertex stabilizers. Its Cayley graph with respect to any finite generating set then has a tree-like structure with finite pieces. Given a finite generating set of $G$, the pieces of the tree-like structure are algorithmically constructed.*

*Therefore, $p_c$ is an algebraic number and one can use the algorithm from Theorem 1.1(*iii*) to compute $p_c$, given any finite generating set.*

PROOF. The group $G$ in this case is a fundamental group of a finite graph of finite groups. Recall that the standard generating set, as defined in Section 6.3, consists of generators of the vertex groups and free letters corresponding to the edges outside the spanning tree of the factor graph $T/G$. Assume that $S_1$ is a maximal generating set which contains all elements of the vertex groups and free letters. In Section 6.3, we constructed a tree-like structure for the Cayley graph, with pieces corresponding to vertex groups. This tree of pieces $P_i$ is a starting point for our generalization to an arbitrary generating set.

Denote by $m$ the number of vertex groups in $T/G$. Elements of an arbitrary (finite) generating set $S$ can be represented by reduced words in $S_1$ of bounded length (depending on $S$). Denote the maximal length by $N$.

The next step is to apply Lemma 2.6 several times so that the piece $P'_i$ contains all of its descendants up to $n$th generation. Every application of that lemma enlarges the pieces by one generation. We will find $n$ large enough to guarantee that there exists a piece containing both endpoints of each edge labeled by an element in $S$.

Let $P'_i$ be a piece in the $n$-times "enlarged" tree-like structure. Then, $V(P'_i) = \bigcup_{k \in \Lambda_i(n)} V(P_k)$, where $\Lambda_i(0) = \{i\}$ and $\Lambda_i(n+1) = \{i\} \cup \bigcup_{s \in \Lambda_i(n)} \Lambda_s$. Since the original pieces were finite, the $\Lambda_j$'s are finite and so are the modified pieces $P'_i$. An intuitive picture of how such enlarged pieces look is provided by the free groups; see Corollary 7.6.

Recall that the pieces $P_i$ correspond to transitions of vertex groups $gG_v$ and that they are connected in the tree of pieces if they are so-called neighbors (corresponding to an edge in the spanning tree) or delayed neighbors (corresponding to other edges in the graph of groups); see Section 6.3. We can now use a graph distance between pieces. Let $s \in S_1$ and $x$ be a vertex in some piece $g_1 G_{v_1}$. How far is a piece containing the vertex $xs$?

Assume $s \in G_{v_2}$ so that it is in some vertex stabilizer. There is a path between $v_1$ and $v_2$ in the spanning tree of the graph $T/G$, visiting vertices $v_1 = u_1, u_2, \ldots, u_k = v_2$. The size of the spanning tree is equal to the the number of vertex groups minus 1, that is, $m - 1$. For every $i$, there is some $g_i \in G$ such that $x \in g_i G_{u_i}$. Therefore, there is a path in the tree of pieces through vertices $g_1 G_{u_1}, g_2 G_{u_2}, \ldots, g_k G_{u_k}$. Moreover, $xs \in g_k G_{v_2}$. Thus, the distance between pieces containing $x$ and $xs$ is at most $m - 1$.

Assume that $s$ is a free letter. It then corresponds to an edge $e$ in the graph of groups outside the spanning tree, starting at vertex $v_2$ and terminating at vertex $v_3$. As before, there is a path from $g_1 G_{v_1}$ to a piece $g_k G_{v_2}$, and $xs \in g_{k+1} G_{v_3}$, a delayed neighbor of $g_k G_{v_2}$. Thus, the distance between pieces containing $x$ and $xs$ is at most $m$, in this case, and similarly for $s^{-1}$.

Let $x \in G$ with word length (in $S_1$) at most $N$. The distance between pieces containing $g$ and $gx$ is at most $mN$ for any $g \in G$. If the distance between two pieces is at most $mN$, then their closest common ancestor differs from them by less than $mN$ generations. In particular, if $\mathrm{dist}(P_i, P_j) \leq mN$, then there exists some $P_k$ such that $P_k'$ contains vertices of $P_i$ and $P_j$, provided $n > mN$.

We can now apply Lemma 2.5 to conclude that the Cayley graph with respect to $S \cup S_1$ has a tree-like structure with finite pieces. Consequently, by Lemma 2.4, the Cayley graph with respect to $S$ also has a tree-like structure with finite pieces. $\square$

**7. Transitive graphs with more than one end.** In what follows, we will generalize the above decomposition into the tree-like structure, which was natural for amalgamated products and HNN extensions, to the case of transitive graphs with more than one end.

Recall that the number of ends of a graph is the supremum of the number of connected components of any of its subgraphs that was obtained by removing a finite set of vertices. An infinite transitive graph can have one, two or infinitely many ends.

7.1. *Dunwoody's result.* First, let us recall the notation and a result obtained by Dunwoody in [3].

Let $\mathcal{G}$ be an infinite connected graph with more than one end. Let $c$ be a subset of vertices of $\mathcal{G}$, let $c^* = V(\mathcal{G}) \setminus c$ and denote by $\partial_E c$ the set of edges having one endpoint in $c$ and the other in $c^*$. Denote by $\partial_V c$ the set of vertices in $c$ having a neighbor in $c^*$. A set of vertices $c$ such that $\partial_E c$ is finite is called a *cut*. A cut is said to be *nontrivial* if both $c$ and $c^*$ are infinite.

LEMMA 7.1 (Dunwoody). *Let $\mathcal{G}$ be a graph with more than one end and let $H \subset Aut(\mathcal{G})$. There exists a nontrivial cut $d$ such that for any $g \in H$, one of the inclusions $d \subset gd$, $d \subset gd^*$, $d^* \subset gd$, $d^* \subset gd^*$ holds.*

Let $d$ be a cut satisfying the above lemma. For any other cut $b$, we define

$$T_b = \{gd \mid g \in Aut(\mathcal{G}), gd \subsetneq b \text{ and there is no } h \in Aut(\mathcal{G})$$
$$\text{such that } gd \subsetneq hd \subsetneq b \text{ or } gd \subsetneq hd^* \subsetneq b\},$$

(16)
$$T_b^* = \{gd^* \mid g \in Aut(\mathcal{G}), gd^* \subsetneq b \text{ and there is no } h \in Aut(\mathcal{G})$$
$$\text{such that } gd^* \subsetneq hd \subsetneq b \text{ or } gd^* \subsetneq hd^* \subsetneq b\},$$

$$Q_1 = d \cap \bigcap_{c \in T_d \cup T_d^*} (c^* \cup \partial_V c),$$

$$Q_2 = d^* \cap \bigcap_{c \in T_{d^*} \cup T_{d^*}^*} (c^* \cup \partial_V c).$$

If there exists some $g \in Aut(\mathcal{G})$ such that $gd = d^*$, the above definition gives us $T_d = T_d^*$. In order not to consider each cut twice, we set $T_d^* = \varnothing$ in this case.

We will now decompose the graph into pieces isomorphic to $Q_1$ and $Q_2$.

LEMMA 7.2. *Let $c_1, c_2 \in T_d \cup T_d^*$ and $c_1 \neq c_2$. The distance from $c_1 \setminus Q_1$ to $c_2 \setminus Q_1$ in $\mathcal{G}$ is then at least 2. Moreover, $d = Q_1 \mathbin{\dot\cup} \bigcup_{c \in T_d \cup T_d^*} (c \setminus \partial_V c)$ and similarly for $d^*$.*

PROOF. From the property that there is no $h \in Aut(\mathcal{G})$ such that $c_i \subsetneq hd \subsetneq d$ or $c_i \subsetneq hd^* \subsetneq d$ for $i = 1, 2$, it follows that $c_1 \not\subset c_2$ and $c_1^* \not\subset c_2^*$. Since $d^* \subset c_i^*$, we have $c_1^* \not\subset c_2$. Therefore, by Dunwoody's result in Lemma 7.1, $c_1 \subset c_2^*$ and $c_2 \subset c_1^*$. If $\partial_V c_i \subset Q_1$, then the claim would follow.

Assume that $c \in T_d \cup T_d^*$ such that $\partial_V c \not\subset Q_1$. This means that there exists some $c' \in T_d \cup T_d^*$ such that $(c'^* \cup \partial_V C') \cap d$ does not contain $\partial_V c$ and therefore $(c' \setminus \partial_V c') \cap \partial_V c \neq \varnothing$. But, from above, we have that $c' \subset c^*$, which is a contradiction. $\square$

THEOREM 7.3. *Every transitive graph with more than one end admits a tree-like structure with finite border sets.*

PROOF. The pieces we consider are $Q_1$, $Q_2$ and their translates. Assume that $d$ is a cut satisfying Lemma 7.1 and that the origin is in $\partial_V d^*$. Let the root piece $P_0$ be a subgraph of $\mathcal{G}$ induced on the vertices $Q_2 \cup \partial_V d^*$.

Now, with each translate of $d$ (resp., $d^*$), we have a translate of $Q_1$ (resp., $Q_2$). Consider the cuts in $T_{d^*} \cup T_{d^*}^*$ and the translates of $Q_1$ and $Q_2$ corresponding to them. The pieces of the first generation are subgraphs induced by the vertices of these translates of $Q_1$ and $Q_2$ (we include all the edges that are not already in $P_0$). Let the first generation also contain a piece induced by edges in the original $Q_1$.

The whole tree of pieces is constructed inductively. We always consider translates of $Q_1$ and $Q_2$ corresponding to the cuts in $T_b \cup T_b^*$, where $b$ are the cuts used in the previous generation. The piece is a subgraph induced on the vertices of the appropriate translate of $Q_1$ and $Q_2$, and we exclude all edges which are already in some piece of a previous generation. The border set $B_i$ is defined as the intersection of $P_i$ with its parent ($B_0 := \{o\}$).

We can now start verifying the properties of Definition 2.1:

(1) The pieces and border sets were defined so that $B_i$ is a subset of vertex set of $P_i$.
(2) Every edge is in exactly one piece—the first appearing in the construction and containing both endpoints of the edge.
(3) The partial order of the elements follows from the construction, as well as the property that a border set is the intersection of a piece with its parent.
(4) If we remove a border set $B_i$, then the graph falls apart. In particular, $U(P_i)$ is the subgraph induced by vertices of $c'$ and by Lemma 7.2, it intersects the $U(P_i)^c$ only by vertices in $\partial_V c$, which are in $B_i$.
(5) There are two isomorphism classes of pieces (excluding the root), based on whether $P_i$ arises from transition $Q_1$ or $Q_2$. □

7.2. *Application to free groups with nonstandard generators.* The simplest example of a transitive graph with infinitely many ends is a regular tree. The Cayley graph of a free group with respect to free generators is a regular tree and the percolation on it is well understood. The result of this section gives us a simple way to find $p_c$ for any finite generating set.

Assume that $F_n = \langle x_1, \ldots, x_n \rangle$. Let $\mathcal{G}$ be its Cayley graph with respect to the standard (free) generators. Denote by $H$ the subgroup of $Aut(\mathcal{G})$ generated by left translations by elements in $F_n$ and isomorphisms arising from permutations of the generators.

The following lemma shows that the the cut $d$ satisfying the condition in Lemma 7.1 can be found explicitly.

LEMMA 7.4. *Consider the free group $F_n$ with any finite set of generators. Denote by $d$ the set of vertices labeled by words starting with the letter $x_1$ (one element of any free generating set). Then, $d$ is such that for any $g \in H$, one of the inclusions $d \subset gd$, $d \subset gd^*$, $d^* \subset gd$, $d^* \subset gd^*$ holds.*

PROOF. For any finite generating set, $d$ and $d^*$ are infinite with finite boundaries. Indeed, for any $m$, there are only finitely many pairs $(a, b)$, $a \in d$, $b \in d^*$, with word distance less than $m$. Without loss of generality, we assume that the generating set is symmetric on the permutation of generators, so the whole group $H$ acts on the Cayley graph by isomorphisms.

Now, for any $h \in H$, if the origin 1 is contained in $hd$, then there is a word $x$ such that all vertices in $hd^*$ are labeled by reduced words starting with $x$ and $x \in d^*$. If the origin is in $hd^*$, then $hd$ has a similar property. Assume that for a given $h \in H$, such a word is $w$ and the origin is in $hd$ (resp., in $hd^*$). Then, if $w$ starts with the letter $x_1$, then $d^* \subset hd^*$ (resp., $d^* \subset hd$), otherwise $d \subset hd$ (resp., $d \subset hd^*$). $\square$

The following theorem is a special case of Theorem 1.5. We include here a different, more "hands on" proof. We present here the pieces of the tree-like structure of free groups explicitly. They are constructed using a cut set satisfying Lemma 7.1. We could equally well use Lemma 2.6 (starting with pieces that are single edges of the Cayley graph of a free group with respect to free generators)— the pieces obtained would be exactly the same.

THEOREM 7.5. *The Cayley graph of free group $F_n$ (with n free generators) with respect to any finite generating set has a tree-like structure with finite pieces. Therefore, its $p_c$ is an algebraic number and one can use the algorithm from Theorem 1.1(*iii*) to compute $p_c$, given any finite generating set.*

PROOF. Recall that we want to show that the Cayley graph of $F_n$ with respect to any finite set of generators has a tree-like structure with finite pieces.

Any finite set of generators is contained in some ball $S_k$, that is, the set of all words of length at most $k$ (in the standard word metric). Therefore, using Lemma 2.4, it is enough to find the tree-like structure with finite pieces for the generating set $S_k$. The isomorphisms will be from the group $H$, which admits the group $F_n$ as a finite-index subgroup.

Consider the cut $d$ from Lemma 7.4. We observed that it is always the case that either $hd$ or $hd^*$ is a set of all vertices labeled by a reduced word starting with some $x$. We denote such a cut by $D_x$ and its complement (containing the origin) by $D_x^*$. For example, for $d$, we write $D_{x_1}$.

Applying the construction in the proof of Lemma 7.3 to the cut $d$, we obtain the following tree-like structure. The pieces are isomorphic to $Q_1$ and $Q_2$. Since $d$ and $d^*$ are isomorphic, it is enough to describe just one of them, say, $Q_2$ [general formula (16)],

$$T_{d^*} \cup T_{d^*}^* = \{D_x : |x| = 1, x \neq x_1\},$$
$$Q_2 = D_{x_1}^* \cap \bigcap_{|x|=1, x \neq x_1} (D_x^* \cup \partial_V D_x).$$

If a vertex is labeled by a reduced word longer than $k$, then it is in $D_x \setminus \partial_V D_x$ for some $|x| = 1$, thus not in $Q_2$. Therefore, $Q_2 \subset S_k$, hence it is finite.

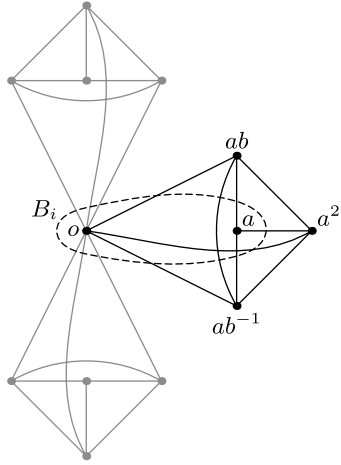

Fig. 4. *Piece in Cayley graph of $F_2$ with respect to $S_2$.*

Note that since $|T_{d^*} \cup T^*_{d^*}| = 2n - 1$, the branching number of the tree of pieces is $2n - 1$. □

The tree-like structure is explicitly given by the procedure of Lemma 7.3. For the Cayley graph of $F_2$ with respect to $S_2$, it results in pieces as in Figure 4. In the picture, we see that the whole $Q_2$ splits into three parts. For the tree-like structure, we can take pieces corresponding to these thirds of $Q_2$. In particular, the middle piece $P_i$ (with labeled vertices as in the picture) has border set $B_i = \{o, a\}$ and its children pieces (three of them) share with $P_i$ the middle vertex $a$ and one of the remaining three vertices $(ab, a^2, ab^{-1})$. The value of the critical probability is $0.139\ldots$ in this case.

If the generating set is a general ball $S_k$, then we can identify the following tree-like structure.

COROLLARY 7.6. *Consider the free group $F_2 = \langle a, b \rangle$ and its Cayley graph with respect to the generating set $S$ containing all words of length at most $k$ (in the standard word metric). There is a tree-like structure with border sets of size $\frac{3^{k-1}+1}{2}$, pieces of size $\frac{3^k+1}{2}$, and the tree of pieces has branching number 3.*

PROOF. The root piece $P_0$ is a subgraph induced by vertices of $S_{k-1}$ with border set containing the origin. The four pieces of first generation are subsets of $S_k$ such that they contain the origin and vertices labeled by words starting with a specific letter. The children of the piece corresponding to the letter $x$ are induced by subsets of $xS_k$ containing $x$ and vertices labeled by words starting with $xy$, where $x \neq y \in \{a, a^{-1}, b, b^{-1}\}$. Clearly, each such piece has exactly three children. Since $|S_k| = 4*(3^k-1)/(3-1)+1$, we can see that $|P_i| = (|S_k|-1)/4+1 = (3^k+1)/2$ for all $i \neq 0$. Again, we included every edge in exactly one piece, that is, in the oldest one containing both of its endpoints. The border set is defined as the intersection of a piece with its parent. It is not difficult to see that $B_i$ coincides with vertices of a piece in the Cayley graph with respect to $S_{k-1}$ and thus $|B_i| = (3^{k-1}+1)/2$. $\square$

**Acknowledgments.** The author would like to thank Mark Sapir, Tatiana Smirnova-Nagnibeda and Russell Lyons for fruitful conversations and useful advice. I am grateful to Vadim Kaimanovich for bringing Dunwoody's result to my attention.

MATHEMATICS DEPARTMENT
VANDERBILT UNIVERSITY
1326 STEVENSON CENTER
NASHVILLE, TENNESSEE 37240
USA
E-MAIL: iva.kozakova@vanderbilt.edu